\documentclass[twoside,12pt,leqno]{article}
\advance\topmargin by -\headheight\advance\topmargin by -\headsep

\usepackage{amsmath}
\usepackage{amsthm}
\usepackage{amssymb}
\usepackage{amscd}
\usepackage{graphicx}
\usepackage{epsfig}

\numberwithin{equation}{section}
\allowdisplaybreaks

\newtheorem{theorem}{Theorem}[section]
\newtheorem{lemma}{Lemma}[section]
\newtheorem{proposition}{Proposition}[section]

\theoremstyle{definition}

\newtheorem{example}{Example} [section]

\newtheorem{remark}{Remark}[section]
\newtheorem*{remarku}{Remark}
\begin{document}
\newcommand{\pp}[2]{\ensuremath{\frac{\partial #1}{\partial #2}}}
\newcommand{\cinfty}{\ensuremath{\mathcal{C}^\infty}}
\newcommand{\pinfty}{\ensuremath{\mathcal{P}^\infty}}
\newcommand{\piv}{\ensuremath{\pi_V}}
\newcommand{\pih}{\ensuremath{\pi_H}}
\newcommand{\Lie}{\ensuremath{\mathcal{L}}}
\renewcommand{\theenumi}{\arabic{enumi}}



\setcounter{page}{1}
\thispagestyle{empty}



\markboth{Benjamin Lent Davis
 and A\"{\i}ssa Wade }{Nonlinearizability
 of certain Poisson structures near a symplectic leaf}

\label{firstpage}
$ $
\bigskip

\bigskip

\centerline{{\Large 
Nonlinearizability
 of certain Poisson structures near a symplectic leaf
  }}

\bigskip
\bigskip
\centerline{\large by Benjamin Lent Davis 
\footnote[1]{Research partially supported by  the Faculty Development Fund of Saint Mary's College.}
and  A\"{\i}ssa Wade
\footnote[2]{Research partially supported by the Shapiro Funds}
 }

\vspace*{.7cm}

\begin{abstract}
We give an intrinsic proof that Vorobjev's first approximation of a
 Poisson manifold near a symplectic leaf is a Poisson manifold.
We also show that Conn's linearization results cannot
 be extended in Vorobjev's setting
\end{abstract}

\pagestyle{myheadings}
\section{Introduction}

A  Poisson structure on a smooth $n$-dimensional manifold $M$ is a bivector
 field $\pi \in \Gamma(\Lambda ^2 TM)$  such that the Schouten bracket
 $[\pi, \pi]=0$.It is known that every Poisson structure $\pi$ on $M$ 
 gives rise to a foliation by symplectic leaves. The study of the local structure of a Poisson manifold near a zero-dimensional leaf 
(i.e., a point  $m$ for which $\pi_m=0$ ) leads to the linearization problem,
 which we now describe. Recall that the {\em linear approximation}
 of a Poisson structure $\pi$ at a zero-dimensional leaf $m$ is determined
 by a Lie algebra $\mathfrak g$, called the {\em  transverse Lie algebra or isotropy Lie algebra} at $m$.
More precisely,  the linear approximation coincides with the 
canonical linear Poisson structure $\pi^{(1)}$ defined on the dual $\mathfrak{g}^*$ of $\mathfrak g$.
In \cite{weinstein83}, Weinstein  showed  that, if $\pi$ is a Poisson structure  on $M$ which vanishes at a point $m$ and whose transverse Lie algebra at $m$ is semi-simple, then $\pi$
 is formally isomorphic to its linear approximation 
$\pi^{(1)}$. The analytic and the smooth versions of this result were
 proved by Conn in \cite{conn84} and \cite{conn85}. Other
 partial results on the linearization of Poisson structures can be found in
 \cite{molinier93} and  \cite{dufour-zung}.  A survey of the literature on linearization of Poisson brackets may be found in \cite{fm04}.

Less studied is  the local structure of a Poisson
 manifold near a symplectic leaf with a nonzero dimension.
 Results on the Poisson topology of neighborhoods of symplectic leaves have been obtained by Ginzburg and Golubev \cite{gg}, Crainic \cite{crainic}, and
Fernandes \cite{fernandes}, \cite{fern02}.
The study of linearization near symplectic leaves of nonzero dimension was initiated by Vorobjev  (see \cite{vorobjev}),
 who defined the anologue of the  linear approximation $(\mathfrak g^*, \pi^{(1)})$
  for a symplectic leaf of dimension
 $d>0$. His goal was to generalize Conn's local linearization result to nonzero dimensional symplectic leaves.

In the present work, we give complete proofs  
of some results (see Theorems \ref{Main 1}
and  \ref{Lie alg}  below)
which were stated in \cite{vorobjev} without proofs.  We also give examples 
 of linearizable and nonlinearizable Poisson structures near a symplectic leaf
  having a nonzero dimension.  In particular, we show that Conn's linearization result cannot be extended in a straightforward way to the semi-local context.  To our knowledge, Example 3 (see below) is the first counter-example to
 the semi-local linearization question.


\section{Integrable geometric data}
Let $S$ be an embedded  submanifold of a smooth $n$-dimensional manifold $P$. Fix a tubular neighborhood  ${\cal N}$ of  $S$.
 This corresponds to a  vector bundle $p : {\cal N} \rightarrow S$. We 
 identify $S$ with the zero section of ${\cal N}$. Moreover, we denote $Vert=\ker p_*$. 
 An \emph{Ehresmann connection} on ${\cal N}$ is a projection map
 $\Gamma: T {\cal N} \rightarrow Vert$. Equivalently, we have 
 a smooth vector subbundle $Hor \subset {\cal N}$ such that
$$T_x{\cal N}= Hor_x \oplus Vert_x \quad \forall \ x \in {\cal N}.$$
 For any vector field $X \in \chi(S)$, there is a unique horizontal
 vector field $\overline{X} \in \chi(Hor)$ which is called the \emph{horizontal lift}
 of $X$, and  satisfies $$ p_*(\overline{X})=X.$$ 
 Define the curvature of $\Gamma$ by the formula
$$Curv_{\Gamma}(X,Y)=[\overline{X}, \overline{Y}]- \overline{[X,Y]}, \quad \mbox{ for any } \ X, Y \in \chi(S) .$$
\noindent  Suppose that $S$ is equipped with a symplectic form $\omega \in \Omega^2(S)$.
  Now, we consider a triple $(\Gamma, \nu, \varphi)$ called  {\em geometric data}
and formed by
 
\smallskip
$\bullet$ an Ehresmann connection $\Gamma$,

\smallskip
$\bullet$  a vertical bivector field $\nu \in \Gamma(\Lambda^2 Vert)$,

\smallskip
$\bullet$  and a nondegenerate 2-form $\varphi \in \Omega^2(S) \otimes \cinfty({\cal N})$ given by
 $$\varphi=\omega \otimes 1 + R,$$
 \noindent where  $\omega$ is the symplectic form on the zero section $S$, and $R$ vanishes on $S$, i.e.
$$\varphi_x(X,Y)=\omega_x(X,Y), \quad \forall \ x \in S, \quad \forall X, Y \in \chi(S).$$

The projection $p:{\cal N}\rightarrow S$ induces a map $p^*:\Omega^2(S)\otimes \cinfty({\cal N})\rightarrow \Omega^2({\cal N})$  by pulling back the first factor of the tensor product.
Define the \emph{coupling 2-form} $\widehat{\omega} \in \Omega^2({\cal N})$ by $\widehat{\omega}=p^*\varphi.
$
\newcommand{\ann}{\ensuremath{\mathrm{ann}}}
In a sufficiently small tubular neighborhood we have that, for any $x \in {\cal N}$, 
the linear map $\widehat{\omega}_{|Hor}(x): Hor_x \rightarrow \ann(Vert_x)$
is an isomorphism, where $\ann(Vert_x)$ is the annihilator of $Vert_x$.. If necessary, we may work with a smaller tubular neighborhood.
Define the \emph{horizontal coupling bivector} by  
\begin{equation}
\label{hcb}
\mu=\widehat{\omega}_{|Hor}^{-1}.
\end{equation}
 This gives  the identity
$$\widehat{\omega}(\mu\alpha,\mu\beta)=-\mu(\alpha,\beta).$$
\noindent With respect to the above notations
 we have the following theorem:
\begin{theorem}
\label{Main 1}
The  bivector field $\pi=\mu +\nu$ is a Poisson bivector field if and only if the 
 following \emph{integrability conditions} are satisfied
\begin{enumerate}
\item $\nu$ is a Poisson  bivector field (i.e. $[\nu, \nu]=0$);
 \item $[\overline{X}, \ \nu]=0,$ for all $X \in \chi(S)$;
 \item $Curv_{\Gamma}(X,Y)= \nu(d\widehat{\omega}(\overline{X}, \overline{Y}))$;
 \item $d \widehat{\omega}(\overline{X}_1, \overline{X}_2, \overline{X}_3)=0$, for all
 $X_1, X_2, X_3 \in \chi(S)$.
\end{enumerate}
\end{theorem}
 
\noindent The proof relies on the following three lemmas.
\begin{lemma}
\label{formula 1}
Let $M$ be a smooth manifold and let $\Lambda$ be a bivector field  on $M$.
 Then, 
$$ -\frac{1}{2}[\Lambda,\Lambda](\alpha,\beta,\gamma)=
\Big(\Lambda(d\Lambda(\alpha,\beta),\gamma)+\langle\alpha,[\Lambda \beta,\Lambda\gamma]\rangle
 \Big)+\mbox{c.p.},$$
\noindent for any $\alpha,\beta,\gamma\in\Omega^1(M))$. 
\end{lemma}
 
\begin{proof}
 We have the following formula (see~\cite{ks-m}):
$$\Lambda({\cal L}_{\Lambda \alpha}\beta - {\cal L}_{\Lambda \beta}\alpha 
-d(\Lambda(\alpha, \beta)))
 = [\Lambda \alpha, \ \Lambda \beta] + \frac{1}{2}[\Lambda,\Lambda](\alpha,\beta, \cdot).$$
\noindent There follows
$$
\frac{1}{2}[\Lambda,\Lambda](\alpha,\beta, \gamma)  =
 - {\cal L}_{\Lambda \alpha} \langle \beta, \ \Lambda \gamma \rangle 
- \langle \alpha, \ [\Lambda \beta, \ \Lambda \gamma] \rangle + \mbox{c.p.}
$$
\end{proof}

\begin{lemma}\label{jacobi} 
 For any horizontal 1-forms $\alpha,\beta,$ and $\gamma$, we have
$$ -\frac{1}{2}[\mu,\mu](\alpha,\beta,\gamma)= d\widehat{\omega}(\mu\alpha,\mu\beta,\mu\gamma).$$
\end{lemma}

\begin{proof} Indeed,
\begin{eqnarray*}
d\widehat{\omega}(\mu\alpha,\mu\beta,\mu\gamma)&=&
\Lie_{\mu\alpha}\Big(\widehat{\omega}(\mu\beta,\mu\gamma)\Big)+\widehat{\omega}(\mu\alpha,[\mu\beta,\mu\gamma])
+\mbox{c.p.}\\
&=&\mu(d\mu(\beta,\gamma),\alpha)+\langle\alpha,[\mu\beta,\mu\gamma]\rangle+\mbox{c.p.}\\
&=& -\frac{1}{2}[\mu,\mu](\alpha,\beta,\gamma)\mbox{ by  Lemma \ref{formula 1}.}
\end{eqnarray*}
\noindent There follows the lemma.
\end{proof}

\begin{lemma}
\label{prep}
For any horizontal 1-form $\alpha$ and for any vertical
 1-forms $\beta, \gamma$, we have
$$-\frac{1}{2}[\pi,\pi](\alpha,\beta,\gamma)= [\overline{X}, \nu](\alpha,\beta),$$
\noindent where $\overline{X}= \mu \alpha$. 
\end{lemma}

\begin{proof}
By Lemma~\ref{formula 1}, we have that
\begin{eqnarray*}
 - \frac{1}{2}[\pi,\pi](\alpha,\beta,\gamma)&=&
\pi\Big(d(\pi(\alpha,\beta)),\gamma\Big)+\pi\Big(d(\pi(\beta,\gamma)),\alpha \Big)
+\pi\Big(d(\pi(\gamma,\alpha)),\beta \Big)\\
&&+\langle
\alpha,[\pi\beta,\pi\gamma]\rangle+\langle\beta,[\pi\gamma,\pi\alpha]\rangle+\langle\gamma,[\pi\alpha,\pi\beta]\rangle.\\
&=& \pi(d\pi(\beta,\gamma),\alpha)+\langle \beta,[\pi\gamma,\pi\alpha]\rangle+\langle \gamma,[\pi\alpha,\pi\beta]\rangle\\
&=& -\mathcal{L}_{\overline{X}}(\pi(\beta,\gamma))+\langle \beta,[\pi\gamma,\pi\alpha]\rangle+\langle \gamma,[\pi\alpha,\pi\beta]\rangle\,\, (\star)
\end{eqnarray*}
since $\pi(\alpha,\beta)=0$ and $\pi(\gamma,\alpha)=0$, and the Lie
 bracket of  two vertical vector fields is again vertical implying that 
$\langle \alpha,[\pi\beta,\pi\gamma]\rangle=0$.
On the other hand, by the product rule,
$$\mathcal{L}_{\overline{X}}(\pi(\beta,\gamma))=(\mathcal{L}_{\overline{X}}\pi)(\beta,\gamma)+\pi(\mathcal{L}_{\overline{X}}\beta,\gamma)+\pi(\beta,\mathcal{L}_{\overline{X}}\gamma)$$
and so
$$(\mathcal{L}_{\overline{X}}\pi)(\beta,\gamma)=\mathcal{L}_{\overline{X}}(\pi(\beta,\gamma))-\pi(\mathcal{L}_{\overline{X}}\beta,\gamma)-\pi(\beta,\mathcal{L}_{\overline{X}}\gamma).\, (\star\star)$$
Concentrating now on the second term of $(\star\star)$, we set $Z=\pi\gamma$ and compute using the product rule that
\begin{eqnarray*}
-\pi(\mathcal{L}_{\overline{X}}\beta,\gamma)&=&(\mathcal{L}_{\overline{X}}\beta)(Z)\\
&=&\overline{X} \cdot \pi(\gamma,\beta)-\beta([\pi\alpha,\pi\gamma])\\
&=&-\mathcal{L}_{\overline{X}}(\pi(\beta,\gamma))+\beta([\pi\gamma,\pi\alpha]).
\end{eqnarray*}   
Performing a similar computation on the third term of $(\star\star)$, we find that
$$(\mathcal{L}_{\overline{X}}\pi)(\beta,\gamma)=-\mathcal{L}_{\overline{X}}(\pi(\beta,\gamma))+\langle\beta,[\pi\gamma,\pi\alpha]\rangle+\langle\gamma,[\pi\alpha,\pi\beta]\rangle.\,\,(\star\star\star)$$

Combining ($\star$) and ($\star\star\star$), we see that 
\begin{displaymath}
\begin{array}{rcl}
- \frac{1}{2}[\pi,\pi](\alpha,\beta,\gamma) &=& (\mathcal{L}_{\overline{X}}\pi)(\beta,\gamma),\\
&=& (\mathcal{L}_{\overline{X}}(\mu+\nu))(\beta,\gamma,)\\
&=& (\mathcal{L}_{\overline{X}}\mu)(\beta,\gamma),\\
&=& 0,
\end{array}
\end{displaymath}
since $\beta$ and $\gamma$ are vertical.
\end{proof}

\medskip

\noindent{\bf Proof of Theorem \ref{Main 1}:} 
The splitting $T{\cal N}=Hor\oplus Vert$ induces the splittings
$$\bigwedge^k T{\cal N}=\bigoplus_{i+j=k} Hor^i \wedge Vert^j,$$
where $Hor^i \wedge Vert^j$ is the wedge product of $i$ copies of the vector bundle $Hor$ with
$j$ copies 
of the vector bundle $Vert$.  A section of $Hor^i \wedge Vert^j$ is said to be 
a \emph{multivector field of degree $(i,j)$}. 

 Thus, the trivector field $[\pi,\pi]$ is a sum of component trivector fields of degrees (3,0), (2,1), (1,2), and (0,3).   We will show that   conditions 1-4
  given in Theorem \ref{Main 1} are
  equivalent to the vanishing of 
 $[\pi,\pi]$ in degree $(n, 3-n)$, where $n=0,..,3$.
We have $$[\pi,\pi]=[\mu+\nu,\mu+\nu] = [\mu,\mu]+2[\mu,\nu]+[\nu,\nu].$$
  In general, the horizontal  distribution  is not necessarily integrable, but 
$Vert$ is always integrable. 

\noindent \textbf{Degree (0,3):} 
\noindent Thus, the degree (0,3) component of $[\pi,\pi]$ is exactly $[\nu,\nu]$. 
\medskip

\noindent \textbf{Degree (1,2):}
\noindent Let $\alpha$ be a horizontal 1-form, and
 let $\beta, \gamma$ be two vertical 1-forms.
By lemma \ref{prep}, we get
\begin{equation}
\label{Eq 2} 
[\pi, \pi](\alpha,\beta,\gamma)=-2 [\overline{X}, \nu](\alpha,\beta),
\end{equation}
\noindent where $\overline{X}= \mu \alpha$. 

\medskip
\noindent \textbf{Degree (2,1):}

\noindent Suppose $\alpha,\beta$ are horizontal 1-forms and $\gamma$ is a vertical 1-form. 
By Lemma~\ref{formula 1}, we have that
\begin{eqnarray*}
- \frac{1}{2}[\pi,\pi](\alpha,\beta,\gamma)&=&
\pi\Big(d\pi(\alpha,\beta),\gamma\Big)+\pi\Big(d\pi(\beta,\gamma),\alpha \Big)
+\pi\Big(d\pi(\gamma,\alpha),\beta \Big)\\
&&+\langle \alpha,[\pi\beta,\pi\gamma]\rangle+\langle\beta,[\pi\gamma,\pi\alpha]\rangle+\langle\gamma,[\pi\alpha,\pi\beta]\rangle.
\end{eqnarray*}
Several terms vanish automatically.  In particular $\pi(\beta,\gamma)=0$ and $\pi(\gamma,\alpha)=0$ since $\pi$ vanishes in degree (1,1). 
 Without loss of generality, we can suppose that $\overline{X}=\mu \alpha$, $\overline{Y}=\mu \beta$
 are horizontal lifts of vector fields $X, Y \in \chi(M)$. Setting
$Z:=\pi\gamma$, we get
$$\langle \alpha,[\pi\beta,\pi\gamma]\rangle=\langle\alpha,\Lie_{\overline{Y}} Z\rangle=0$$
and 
$$\langle \beta,[\pi\gamma,\pi\alpha]\rangle=\langle\beta,-\Lie_{\overline{X}} Z\rangle=0$$
since the Lie derivative of a vertical vector field by the horizontal lift of a vector field on $S$ is again vertical.
Thus,
\begin{eqnarray*}
- \frac{1}{2}[\pi,\pi](\alpha,\beta,\gamma)
&=&\pi(d(\pi(\alpha,\beta)),\gamma)+\langle\gamma,[\pi\alpha,\pi\beta]\rangle\cr
&=& \nu(d(\mu(\alpha, \beta)), \gamma))+\langle\gamma, [\overline{X}, \overline{Y}]\rangle\cr
&=&-\nu(d(\widehat{\omega}(\overline{X}, \overline{Y})), \gamma)
+\langle\gamma, [\overline{X}, \overline{Y}]\rangle
\end{eqnarray*}
\noindent But
$$[\overline{X}, \overline{Y}]=\overline{[X,Y]} +Curv_{\Gamma}(X,Y) \quad
{\rm and} \quad \langle \gamma, \overline{[X,Y]}\rangle=0.$$
\noindent It follows
\begin{equation}
\label{Eq 3}
- \frac{1}{2}[\pi,\pi](\alpha,\beta,\gamma)= -
\nu(d(\widehat{\omega}(\overline{X}, \overline{Y})), \gamma)+ 
\langle \gamma, Curv_{\Gamma}(X,Y)\rangle.
\end{equation}

\noindent 
\textbf{Degree (3,0):} 
\noindent We have seen earlier that the degree (3,0) components of $[\pi,\pi]$ and $[\mu,\mu]$ are equal, thus
$$[\pi,\pi](\alpha,\beta,\gamma)=[\mu,\mu](\alpha,\beta,\gamma).$$ 
By Lemma~\ref{formula 1}, 
\begin{equation}
\label{Eq 4}
d\widehat{\omega}(\mu\alpha,\mu\beta,\mu\gamma)=-\frac{1}{2}[\mu,\mu](\alpha,\beta,\gamma).
\end{equation}

\noindent This completes the proof of Theorem \ref{Main 1}.
\hfill $\Box$

\medskip
\noindent{\bf Definition} A geometric data $(\Gamma, \nu, \varphi)$ satisfying conditions 1-4 of Theorem \ref{Main 1}
 is said to be  {\em integrable}. In  this case, the corresponding  bivector field $\pi=\mu + \nu$ 
  is called the {\em  coupling Poisson bivector field associated with} $(\Gamma, \nu, \varphi)$. 
\medskip

\noindent{\bf Remark 1}
Let $X_1, X_2 , X_3$ be vector fields on the base manifold $S$,
 and let $\overline{X}_1, \overline{X}_2, \overline{X}_3$ be their  horizontal lifts respectievely.
 Using the above notations, we have
$$
d\widehat{\omega}(\overline{X}_1, \dots, \overline{X}_3)=
 \Lie_{\overline{X}_1}\Big(\widehat{\omega}(\overline{X}_2, \overline{X}_3)\Big)-
\widehat{\omega}([\overline{X}_1, \overline{X}_2], \overline{X}_3)+c.p.  $$
\noindent The curvature $Curv_{\Gamma}$ takes values in the vertical bundle $Vert$.
It follows that
$$\widehat{\omega}\Big([\overline{X}_1, \overline{X}_2], \overline{X}_3\Big)=
\widehat{\omega}\Big(\overline{[X_1, X_2]}, \overline{X}_3\Big).$$
\noindent Consequently, 
\begin{eqnarray*}
d\widehat{\omega}(\overline{X}_1, \dots, \overline{X}_3)&=&
\Lie_{\overline{X}_1}\Big(\widehat{\omega}(\overline{X}_2, \overline{X}_3)\Big)-
\widehat{\omega}\Big(\overline{[X_1, X_2]}, \overline{X}_3\Big)+c.p. \cr
&=& \Lie_{\overline{X}_1}\Big(\varphi(X_2, X_3)\Big)-\varphi\Big([X_1, X_2], X_3\Big)+c.p.\cr 
&=& (\partial_{\Gamma}\varphi)(X_1,X_2, X_3),
\end{eqnarray*}

\noindent where $$\partial_{\Gamma}:\Omega^k(S)\otimes\cinfty({\cal N})
\rightarrow\Omega^{k+1}(S)\otimes\cinfty({\cal N})$$ is the operator defined by the formula
\begin{eqnarray*}
\partial_{\Gamma}\mathbb F(X_0,\ldots,X_k)&=&\sum_{i=0}^n (-1)^i
 \Lie_{\overline{X}_i}\Big(\mathbb F(X_0,\ldots,\hat{X_i},\ldots,X_k)\Big)\\
&&+\sum_{i<j} (-1)^{i+j}\mathbb F([X_i,X_j],X_0,\ldots,\hat{X_i},\ldots,\hat{X_j},\ldots,X_k).\end{eqnarray*}
\noindent Hence integrability condition 4 in Theorem \ref{Main 1} can be replaced by $\partial_{\Gamma}\varphi=0$.

\medskip
\noindent{\bf Remark 2}
We have shown that, given a manifold $P$, any  integable geometric data $(\Gamma, \nu,\varphi)$ with respect
to a tubular  neighborhood ${\cal N}$ of a symplectic submanifold $(S,\omega)$ of $P$
induces a Poisson bivector field $\pi$ on ${\cal N}$,
 which admits $S$ as a symplectic leaf. The converse is also true.
Indeed, suppose $\pi$ is a Poisson structure on the total space
${\cal N}$ of a vector bundle over a symplectic manifold $(S, \omega)$,
 where $S$ is identified with the zero section of ${\cal N}$, and  $S$
 is a symplectic leaf of ${\cal N}$.
 Let $Vert=\ker p_*$, where $p$ is the natural projection of ${\cal N}$ onto $S$. In addition,
 we assume that $\pi$ is {\em horizontally nondegenerate}, that is
$\pi_{|\ann(Vert)}$ is nondegenerate.
 Then, there is a natural Ehresmann connection
 associated with $\pi$ which is determined by the horizontal subbundle
 $$Hor= \pi(\ann(Vert)).$$
 \noindent The bivector field $\pi$ can be decomposed into
$\pi=\mu+\nu$, where $\mu$ and $\nu$ are horizontal and
 vertical, respectively. Furthermore, we can define
$$\varphi(X,Y)= \pi^{-1}_{|Hor} (\overline{X}, \overline{Y}),$$
\noindent where $\overline{X}, \overline{Y}$ are the lifts of the vector
 fields $X, Y \in \chi(S)$. Using Equations (1)--(4), we
  conclude that integrability conditions 1--4 are satisfied.

\medskip
\noindent {\bf Remark 3.}
Brahic \cite{brahic} and Vaisman \cite{vaisman04} have each given independent proofs of  Theorem \ref{Main 1}.    The proof given here
 is intrinsic, and  has the advantage of clearly relating Vorobjev's linearization formula to geometric identities for transitive Lie algebroids
 (see below).

\section{Transitive Lie algebroids}

A \emph{Lie algebroid} over a manifold $S$ is a real vector bundle
  $p: E \rightarrow S$ together with a bundle map 
$\rho: E\rightarrow TS$ and a real Lie algebra structure 
$[\cdot,\cdot]_E$ on $\Gamma(E)$ such that the following 
 \emph{Leibniz rule} holds
$$[v,fw]_E=f[v,w]_E+(\rho(v)\cdot f)w,$$
\noindent for any $f\in\cinfty(S)$ and $v,w\in\Gamma(E)$.
The map $\rho$ is called the \emph{anchor} of the Lie algebroid.  
Using the Jacobi identity for $[ \cdot , \cdot]_E$
 and the Leibniz identity, one can show that
 the induced map $\rho:\Gamma(E)\rightarrow\chi(S)$ is a Lie algebra 
 homomorphism, i.e.
$$\rho[u,v]_E=[\rho(u), \rho(v)].$$

\medskip

Given any Lie algebroid, the distribution spanned by the image of the anchor is integrable,
 and the leaves of the resulting foliation are called the \emph{orbits} of the Lie algebroid.
  When the anchor is surjective, we call the Lie algebroid \emph{transitive}.

If $E\stackrel{\rho}{\rightarrow} S$ is a transitive Lie algebroid  with anchor map $\rho$,
then we obtain a short exact sequence of vector bundles
$$0\rightarrow I\stackrel{v}{\rightarrow} E\stackrel{\rho}{\rightarrow} TS\rightarrow 0,$$
where $I=\ker\rho$ is called the \emph{isotropy bundle}.

\medskip
 Recall that a \emph{connection} for a transitive Lie algebroid $E\stackrel{\rho}{\rightarrow} S$ 
 is a  splitting $\sigma:TS\rightarrow E$ of the sequence above.
Its  \emph{curvature } $R_\sigma\in\Omega^2(S)\otimes\Gamma(I)$ is given  by
$$R_\sigma (X,Y)=[\sigma X,\sigma Y]-\sigma [X,Y], \quad \mbox{ for all $X,Y\in\chi(S)$}.$$

Any connection $\sigma$ for a transitive Lie algebroid induces a covariant derivative 
for the isotropy bundle $I$ defined by
$$\nabla^{\sigma}_X s=[\sigma X,s], \quad \mbox{ for all $X\in\chi(S)$ and $s\in\Gamma(I)$.}$$
By a  \emph{covariant derivative} for a vector bundle ${\cal N}\rightarrow S$,
 we mean a linear map $\nabla:\Gamma({\cal N})\rightarrow\Omega^1(S)\otimes\Gamma({\cal N})$ satisfying 
the following properties:
$$\nabla_X fs= X(f)\cdot s+ f\nabla_X s \quad
 {\rm and} \quad \nabla_{fX} s=  f\nabla_X s \,$$ for all $X\in\chi(S)$, $f\in\cinfty(S)$, 
and $s\in\Gamma({\cal N})$.

Let $\sigma$ be a connection for a transitive Lie algebroid  $E$.
For simplicity, the corresponding covariant derivative
 will be denoted by $\nabla$ instead of $\nabla^{\sigma}$ when there is no ambiguity.
Define the \emph{curvature of} $\nabla$  by 
$$R_{\nabla}(X_1,X_2)=[\nabla_{X_1} , \ \nabla_{X_2}] 
-\nabla_{[X_1,X_2]} ,$$
\noindent for all $X_1,X_2\in\chi(S)$. 

The Jacobi identity yields geometric identities when evaluated on isotropic and coisotropic sections.
\begin{theorem}
\label{connection}
 Let $\sigma: TS\rightarrow E$ be a connection for a transitive Lie algebroid 
and let $\nabla$ be the induced covariant derivative for the isotropy bundle $I$. 
Then, for all $X, X_1, X_2, X_3\in\chi(S)$ and for all  $s, s_1, s_2\in\Gamma(I)$, we have
\begin{enumerate}
\item []{\bf (i)} The isotropy bundle $I$ is a Lie algebroid.
\item []{\bf (ii)} $\nabla_X[s_1,s_2]=[\nabla_X s_1,s_2]+[s_1,\nabla_X s_2]$,
\item []{\bf (iii)} $[R_{\sigma}(X_1,X_2),s]-R_{\nabla}(X_1,X_2)(s)=0$,
\item []{\bf (iv)} $\nabla_{X_1}R_{\sigma}(X_2,X_3)+R_{\sigma}(X_1,[X_2,X_3])+\mbox{c.p.}=0.$
\end{enumerate}
\end{theorem}
\begin{proof}
Simply use the definitions of the connection, 
the covariant derivative and their respective curvatures, and the fact that 
$$\mathcal{J}(\sigma X,s_1,s_2)=0, \quad \mathcal{J}(\sigma X_1, \sigma X_2,s)=0,
 \quad{and} \quad \mathcal{J}(\sigma X_1, \sigma X_2, \sigma X_3)=0,$$
 where $\mathcal{J}$ is the jacobiator of the bracket on sections of $E$.
\end{proof}
\noindent Versions of this theorem have been noted by several authors
  (see~\cite{mackenzie},\cite{fernandes},\cite{vorobjev}).\medskip
\medskip

 Any connection $\sigma$ for a transitive Lie algebroid $E \rightarrow S$
 gives rise to a {\em homogeneous} Ehresmann connection $\Gamma$ on the dual $I^*$ of the isotropy
bundle, i.e. the horizontal lift of every vector field $X \in \chi(S)$ preserves
 the space $\cinfty_{\mbox{lin}}({I^*})$ of fiberwise linear functions on $I^*$.
 To define the horizontal lift of $X$, we introduce 
 the natural isomorphism 
\begin{equation}
\label{ell1}
\ell: \Gamma(I)\rightarrow\cinfty_{\mbox{lin}}(I^*)
\end{equation}
given by the natural pairing. The horizontal lift 
 $\overline{X}$ of a vector $X\in \chi(S)$ satisfies
$$ \Lie_{\overline{X}}(\ell(s))=\ell(\nabla_X s).$$
Let $Hor \subset TI^*$ be the corresponding horizontal subbundle. Then
 $$TI^*=Hor \oplus Vert,$$
\noindent where $Vert=\ker  p_*$  with $p_*: TI^* \rightarrow TS$.
Define the \emph{Lie-Poisson bivector} $\nu$  by the following formula:
\begin{equation}
\label{Eq 6}
\nu(d u_1, du_2)=\ell([s_1, s_2]), \quad \mbox{where}  \ \ell(s_i)=u_i.
\end{equation}
The Lie-Poisson bivector defines a Poisson structure on $TI^*$ 
 (see Section 16.5 of \cite{weinstein}).
 Now, we suppose that $S$ is equipped with a symplectic form $\omega$.
 Let 
 \begin{equation}
 \label{Eq 7}
 \varphi= \omega\otimes 1 + \ell\circ R_{\sigma}.
 \end{equation}
\begin{theorem}
\label{Lie alg}
 Given a transitive Lie algebroid $E$ over a symplectic manifold $S$ 
 together with
a connection $\sigma: TS \rightarrow E$,
the geometric data $(\Gamma, \nu, \varphi)$ defined as above is integrable.
\end{theorem}

\begin{proof}
We will show that geometric identities \textbf{(i)--(iv)}
  in Theorem \ref{connection} are equivalent to integrability
  conditions 1--4 in Theorem \ref{Main 1}.

\noindent \textbf{1.}  Immediate by Equation \ref{Eq 6}.

\noindent \textbf{2.} By Theorem \ref{connection}\textbf{(ii)}, we have 
$$\nabla_X[s_1,s_2]-[\nabla_X s_1,s_2]-[s_1,\nabla_X s_2]=0.$$
\noindent But
$$\ell(\nabla_X[s_1,s_2])=\langle \overline{X}, d \ell([s_1, s_2])
\rangle=\Lie_{\overline{X}}(\nu(du_1, du_2)),$$
where $\ell(s_i)=u_i$. Similarly
$$\ell([\nabla_X s_1,s_2]+[s_1,\nabla_X s_2])
=\nu(d( \Lie_{\overline{X}}  u_1), du_2)+
\nu(du_1,d( \Lie_{\overline{X}}  u_2)).$$
\noindent It follows that 
$$\ell(\nabla_X[s_1,s_2]-[\nabla_X s_1,s_2]-[s_1,\nabla_X s_2])=
 [\overline{X}, \nu] (du_1, du_2).$$
\noindent Hence,
$$
\nabla_X[s_1,s_2]-[\nabla_X s_1,s_2]+[s_1,\nabla_X s_2]=0
\iff [\overline{X}, \nu]=0.
$$

\noindent \textbf{3.} By Theorem \ref{connection}\textbf{(iii)}, we have 
\begin{eqnarray*}
0&=&\ell \Big([R_{\sigma}(X_1,X_2), s]-R_{\nabla}(X_1,X_2)(s)\Big)\cr
 &=&\nu\Big(d(\ell(R_{\sigma}(X_1, X_2)), \ d(\ell(s)) \Big)-
 \langle [ \overline{X}_1, \overline{X}_2]+\overline{[X_1, X_2]}, d \ell(s) \rangle.
\end{eqnarray*}
\noindent Since $\varphi(X_1, X_2)- \ell(R_{\sigma}(X,Y))$ is constant on each
 fiber, we have
$$\nu\Big(d(\ell(R_{\sigma}(X_1, X_2)), \ d(\ell(s)) \Big)=
\nu\Big(d(\varphi(X_1, X_2)), \ d(\ell(s)) \Big).$$
 Therefore,
 \begin{eqnarray*}
0&=& \ell \Big([R_{\sigma}(X_1,X_2), s]-R_{\nabla}(X_1,X_2)(s)\Big)\cr
 &=&\nu\Big(d(\varphi(X_1, X_2)), \ d(\ell(s)) \Big)-
 \langle Curv_{\Gamma}(X_1, X_2), \ d \ell(s) \rangle \cr
 &=&\nu\Big(d(\widehat{\omega}(X_1, X_2)), \ d(\ell(s)) \Big)-
 \langle Curv_{\Gamma}(X_1, X_2),  \ d \ell(s) \rangle .
 \end{eqnarray*}
We obtain that
$$
 ad_{R_{\sigma}(X_1,X_2)}-R_{\nabla}(X_1,X_2)=0 \iff
Curv_{\Gamma}(X_1, X_2)=
\nu^\sharp d(\widehat{\omega}(X_1, X_2)). 
$$

\noindent \textbf{4.} By Theorem \ref{connection}\textbf{(iv)}, 
\begin{eqnarray*}
\ell\Big(\nabla_{X_1}R_{\sigma}(X_2,X_3)+R_{\sigma}(X_1,[X_2,X_3])\Big)
&=&  \Lie_{\overline{X_1}}\Big( \ell \circ R_{\sigma}(X_2,X_3)\Big)\cr
&& +\ell \circ R_{\sigma}([X_1,X_2],X_3).\cr
\end{eqnarray*}
\noindent Since $\partial_{\Gamma} (\omega \otimes 1)= d \omega \otimes 1=0$, we get
\begin{eqnarray*}
\ell\Big(\nabla_{X_1}R_{\sigma}(X_2,X_3)+R_{\sigma}(X_1,[X_2,X_3])\Big)+c.p.
&=&\Lie_{\overline{X_1}}\Big( \varphi (X_2,X_3)\Big)\cr
 &&+\varphi([X_1,X_2],X_3) +c.p\cr
& & \cr
 &=& \partial_{\Gamma} \varphi (X_1, X_2, X_3).
\end{eqnarray*}
This completes the proof. \end{proof}

The resulting coupling Poisson structure on $I^*$ depends on the choice of connection, but is unique up to isomorphism by the following proposition of Vorobjev proved in \cite{vorobjev}.

 \begin{proposition}
Let  $E_1$ and $E_2$ be two isomorphic transitive Lie algebroids over the same 
 symplectic base manifold $(S,\omega)$,
 and let $\sigma_1: TS \rightarrow E$, $\sigma_1: TS \rightarrow E$ be 
 two connections.  There exists a diffeomorphism
  $\psi$ from a neighborhood ${\cal V}_1$ of the zero section
 $S \subset I_1^*$ onto a neighborhood ${\cal V}_2$ of the zero section
 $S \subset I_2^*$ such that $\psi_{|S}=id$ and
 $\psi_*\pi_1=\pi_2$, where $I_i$ is the isotropy bundle of $E_i$ and
 $\pi_i$ is the coupling Poisson bivector field associated with $\sigma_i$, for $i=1,2$.
  
\end{proposition}

\section{Vorobjev (non)linearizability}

In this section, we are interested in the particular case where the isotropy bundle is the conormal bundle $N^*S$
 to a symplectic leaf $S$ of a Poisson manifold.
Given a Poisson manifold $(P,\Lambda)$, the Poisson tensor
induces a natural map from the cotangent bundle to the tangent bundle by 
 the formula
\begin{eqnarray*}
\Lambda ^\sharp :T^*P&\rightarrow&TP\\
\alpha&\mapsto& \Lambda(\alpha,\cdot).\\
\end{eqnarray*}

It is known that  $T^*P \rightarrow P$ is a Lie algebroid,
 called the \emph{Poisson algebroid} (see \cite{weinstein},\cite{vaisman})
whose  anchor map is $\Lambda^\sharp$ and whose Lie bracket is given by
$$[\alpha, \beta]_{T^*P}= \Lie_{\Lambda^\sharp \alpha} \beta -
\Lie_{\Lambda^\sharp \beta} \alpha - d\Lambda(\alpha, \beta).$$

The restriction to a symplectic leaf $(S, \omega)$ give a transitive
 Lie algebroid $T^*P_{|S} \rightarrow S$. Furtheremore, the kernel of the anchor map
 $\rho$ of this transitive Lie algebroid
 coincides with the conormal bundle $N^*S$ of $S$. 
Let $E$ be a tubular neighborhood of $S$, and let $p: E \rightarrow S$ be  the corresponding
vector bundle. 
The derivative of $p$ gives a connection for the 
transitive Lie algebroid $T^*P_{|S} \rightarrow S$, namely,
$\sigma:TS\rightarrow T^*P|_S$ defined by
$$\sigma(X)=p^*\omega(X),$$
where $\omega\in\Omega^2(S)$ is symplectic form of $S$.
This connection is called the \emph{pullback connection} induced by the tubular neighborhood.
We have the splitting
$$TP_{|S}=TS \oplus NS.$$ 
\noindent Moreover, we know that the normal bundle $NS$ is endowed with a canonical Poisson
  structure
 $\nu$ (see Equation (5)).  Theorem \ref{Lie alg} says that, up to a shrinking of the tubular
  neighborhood $E$ of $S$,
 there is a coupling Poisson bivector field $\pi$ on $E$ having $S$ as a symplectic leaf.

\medskip
\noindent {\bf Definition.}
 The Vorobjev-Poisson structure $\pi$ is called
 the {\em first approximation} of $\Lambda$ at the symplectic leaf 
$(S, \omega)$ with respect to the neighborhood $E$. 

We say that $\Lambda$ is {\em Vorobjev linearizable at the symplectic leaf $S$} if 
 there is a tubular neighborhood $E \subset P$ of $S$ with fibers tangent to the
normal bundle $NS$ and a diffeomorphism $\psi : E \rightarrow U \subset P$ such
$\psi_{|S}=id$ and $\psi_* \pi= \Lambda$.

We now introduce a family of Poisson manifolds called Casimir-weighted products 
and use them to give some concrete examples of Vorobjev \\
(non)linearizability.

Let $(P_1,\pi_1)$ and $(P_2,\pi_2)$ be Poisson manifolds with Casimir functions $f_1$ and $f_2$, respectively. 
 The \emph{Casimir-weighted product} is the Poisson manifold $(P_1\times P_2,
f_2\pi_1+f_1\pi_2)$. 
 That $f_2\pi_1+f_1\pi_2$ is Poisson is an easy computation using Schouten brackets. 
 The next proposition describes the symplectic leaves of a Casimir-weighted product.
\begin{proposition}
\label{productleaves}
Suppose that $(P_1,\pi_1)$ and $(P_2,\pi_2)$ are Poisson manifolds with smooth nowhere
vanishing
 Casimir functions $f_1$ and $f_2$, respectively.  
Then every symplectic leaf of the Casimir-weighted product is of the form
$(S_1\times S_2,\frac{1}{f_2}\omega_1+\frac{1}{f_1}\omega_2)$, where $(S_1,\omega_1)$ and
 $(S_2,\omega_2)$ are symplectic leaves of $(P_1,\pi_1)$ and $(P_2,\pi_2)$, respectively.
\end{proposition}

We now compute Vorobjev linearizations of  Casimir-weighted products of  symplectic Poisson manifolds by  Lie-Poisson manifolds.
\begin{theorem}\label{vorobcasimir}
Let $(S,\pi_S)$ be a symplectic Poisson manifold, and let $f$ be a Casimir for a Lie-Poisson
 manifold $(\mathfrak{g}^*,\pi_{\mathfrak{g}^*})$ such that $f(0)=1$.  Then the Vorobjev linearization 
of the Casimir-weighted product at the leaf $S\times\{0\}$ is
 $$\left(S\times\mathfrak{g}^*, \frac{1}{J^1_0f} \pi_S+\pi_{\mathfrak{g}^*}\right),
$$where $J^1_0f=1 + d_0f \in\cinfty(S\times \mathfrak{g}^*)$ denotes the first jet of $f$ at 0.
\end{theorem}
\begin{proof}
Projection to the first factor makes the Casimir-weighted product into a vector bundle 
$pr_1:S\times\mathfrak{g}^*\rightarrow S$.   Any vector bundle is canonically isomorphic to 
the normal bundle of the zero-section, thus $S\times\mathfrak{g}^*\tilde{\rightarrow} NS$
 and so we have  a canonical tubular neighborhood of the symplectic leaf $S\times {0}$.

We now find the coupling Poisson bivector $\pi=\mu+\nu$ induced by this tubular neighborhood. 
By Equation (\ref{Eq 6}), the Lie-Poisson bivector is the Poisson bivector of the
direct product of $(S,0)$ and the Lie-Poisson manifold $(\mathfrak{g}^*,\pi_{\mathfrak{g^*}})$.
Essentially this means that $\nu=\pi_{\mathfrak{g^*}}$.

By Equation (\ref{hcb}),  the computation of the horizontal coupling form $\mu$ is a two step processes: First we must compute the horizontal distribution $Hor$ of the connection, and then we must compute the coupling 2-form $\widehat{\omega}$.

We now show that the covariant derivative $\nabla$ on the normal bundle $NS$ is 
flat and that the horizontal distribution $Hor$ has leaves $S\times\{\xi\}$ for each
$\xi\in\mathfrak{g}^*$. Let $(q_1,\ldots,q_n,p_1,\ldots,p_n)$ be a Darboux chart on $S$. 
Let  $y_1,\ldots,y_k$ for $\mathfrak{g}$ be any basis.  Each $y_j$ is a linear function on $\mathfrak{g}^*$
 representing a section $dy_j$ of $N^*S$ by pulling-back by the composition of maps 
$S\times\mathfrak{g}^*\rightarrow\mathfrak{g}^*\tilde{\rightarrow}T_0\mathfrak{g}^*$.
Then $$(pr_1^* q_1, \ldots , pr_1^*q_n, pr_1^* q_1, \ldots , pr_1^*q_n, y_1 , \ldots , y_m)$$
 are coordinates on $S \times \mathfrak{g}^*.$
We see that each $dy_j$ is a $\nabla$-parallel section by the computation
\begin{displaymath}
\begin{array}{rcl}
\nabla_{X_{q_i}}dy_j&=&[\sigma(X_{q_i}),dy_j]\\
&=& [pr_1^*dq_i,dy_j]\\
&=& d\{pr_1^*q_i,y_j\}\\
&=& 0,
\end{array}
\end{displaymath}
since $q_i$ and $y_j$ are coordinates on the first and second factor of $S\times\mathfrak{g}^*$,
 respectively, and the bracket $\{\, ,\, \}$ is computed using the bivector $f\pi_S+\pi_{{\mathfrak{g}^*}}$.
Thus, $\{dy_1,\ldots,dy_k\}$ is a $\nabla$-parallel frame field for $N^*S$.
 The dual frame $\{\partial y_1,\ldots,\partial y_k\}$ is necessarily parallel for the dual covariant 
derivative on the dual bundle $NS$.
 Consequently, the horizontal distribution  for $NS$ is spanned by 
$\{\partial q_1,\ldots,\partial q_n,\partial p_1,\ldots,\partial p_n\}$, and so the leaves of the distribution are as claimed above.

Equation (\ref{Eq 7}) for the horizontal coupling form is $\varphi=\omega\otimes 1 + \ell\circ R_\sigma$.  Thus, we must compute the curvature of the pullback connection, $R_\sigma$. 
 In particular, 
\begin{displaymath}
\begin{array}{rcl}
\ell\circ R_\sigma(X_{q_i},X_{p_j})&=&[\sigma(X_{q_i}),\sigma(X_{p_j})]-\sigma([X_{q_i},X_{p_j}])\\
&=& \ell([pr_1^*dq_i,pr_1^*dp_j])\\
&=&\ell(d\{pr_1^*q_i,pr_1^*p_j\})\\
&=&\delta_{ij} \ell(df) \\
&=& \delta_{ij} \Big({\partial f \over \partial y_k} (0) \ell (d y_k) \Big),
\end{array}
\end{displaymath}
showing that $\ell \circ R_\sigma=
\omega\otimes (J^1_0f -1),$ where $J_0^1f$ is the 1-jet of $f$ at $0$.
Consequently, the coupling 2-form is $\varphi=\omega\otimes 1 +\ell \circ R_\sigma=\omega\otimes J^1_0f,$
and so the coupling bivector is
$\mu=\frac{1}{J^1_0f} \pi_S.$
\end{proof}

\newcommand{\pu}{\ensuremath{\partial_u}}
\newcommand{\pv}{\ensuremath{\partial_v}}

\begin{example}{\textbf{1}}
Let $(T,\pi_T=\pu\wedge\pv)$ be the  \emph{unit symplectic torus} obtained by identifying
opposite sides 
of the unit square. Let $(\mathfrak{g}^*,0)=(\mathbb{R},0)$ be the 1-dimensional Lie-Poisson manifold, and let $f(z)=e^z$ be a Casimir.  The Casimir-weighted product 
$$\left(T\times\mathbb{R},e^z\pu\wedge\pv\right),$$
has Vorbjev linearization given by 
$$\left(T\times\mathbb{R},\left(\frac{1}{1+z}\right)\pu\wedge\pv\right),$$
A  linearizing isomorphism 
is given by $\psi(x,y,z)=(x,y,-1+e^{-z}).$
\end{example}

  \begin{example}{\textbf{2}}
Let $(T,\pi_T=\pu\wedge\pv)$ be the unit symplectic torus.
Let $(\mathfrak{g}^*,\pi_\mathfrak{g})=(\mathbb{R},0)$ be the 1-dimensional Lie-Poisson
manifold, and 
let $f(z)=1+z^2$ be a Casimir.  The Casimir-weighted product 
$$\left(T\times\mathbb{R},(1+z^2)\pu\wedge\pv\right),$$
has Vorbjev linearization given by 
$$\left(T\times\mathbb{R},\pu\wedge\pv\right).$$
There is no Poisson isomorphism from the Casimir-weighted product to the Vorobjev linearization.  To see this, note that any such isomorphism must induce isomorphisms of symplectic leaves.  Recall that an invariant of a compact 2n-dimensional symplectic manifold is the \emph{symplectic volume},
$$\mathrm{Vol}(S,\omega)=\int_S \omega^n.$$
The Casimir-weighted product  has leaves of non-constant symplectic volume $1+z^2$, but the Vorobjev linearization has leaves of constant symplectic volume 1.
\end{example}

The following example shows that a Poisson manifold is not necessarily Vorobjev linearizable at a
symplectic leaf possessing 
a semisimple transverse Lie algebra of compact type. In other words, Conn's theorem
 \cite{conn85} cannot be extended in this context.
\newcommand{\px}{\ensuremath{\partial_x}}
\newcommand{\py}{\ensuremath{\partial_y}}
\newcommand{\pz}{\ensuremath{\partial_z}}

\begin{example}{\textbf{3}}
Let $(T,\pi_T=\pu\wedge\pv)$ be the unit symplectic torus. Consider the Lie-Poisson manifold
$(\mathfrak{g}^*,\pi_{\mathfrak{g}^*})=(\mathfrak{so}(3)^*,x\py\wedge\pz+y\pz\wedge\px+z\px\wedge\py).$  Recall that the symplectic leaves of $\mathfrak{so}(3)^*$ are the origin together with the spheres centered at the origin, and so $f=1+x^2+y^2+z^2$ is a Casimir.  Moreover, the sphere at radius $r$ has symplectic form $\omega=\frac{1}{r}dA$, where $dA=dx\wedge dy+dy\wedge dz+dz\wedge dx$ is the standard area form on $\mathbb{R}^3$ (see\cite{mr}, pp.457-458).
By Proposition \ref{productleaves}, $S=T\times\{0\}$ is a symplectic leaf of the Casimir-weighted product
$$(P,\pi)=(T\times\mathfrak{so}(3)^*,f\pi_T+\pi_{\mathfrak{g}^*}).$$
 We will show that $(P,\pi)$ is not Vorobjev linearizable at $S$.
By Theorem \ref{vorobcasimir}, the Vorobjev linearization of $(P,\pi)$ at $S$ is the direct product
$$(T\times\mathfrak{so}(3)^*,\pi_T+\pi_{\mathfrak{g}^*}).$$
Suppose for a contradiction that $$\psi:(T\times\mathfrak{so}(3)^*,f\pi_T+\pi_\mathfrak{g})\tilde{\rightarrow}(T\times\mathfrak{so}(3)^*,\pi_T+\pi_{\mathfrak{g}^*})$$
is an isomorphism of Poisson manifolds.  The map $\psi$ induces isomorphisms of symplectic leaves.   By Proposition \ref{productleaves},
$$(S_1,\omega_1)=(T\times S^2_{r_1},\frac{1}{1+r_1^2}\omega_T+\frac{1}{r_1}dA)$$
is a 4-dimensional symplectic leaf of $(P,\pi)$ when $r_1>0$.
Again by Proposition \ref{productleaves}, the only 4-dimensional leafs of the direct product are of the form
$$(S_2,\omega_2)=(T\times S^2_{r_2},\omega_T+\frac{1}{r_2}dA)$$  
when $r_2>0$.
Since $\psi$ restricts to a symplectomorphism
$(S_1,\omega_1)\tilde{\rightarrow}(S_2,\omega_2),$ there must be leaves of  equal symplectic volume.  The symplectic volume of a product is the product of the symplectic volumes, thus,
$$\mathrm{Vol}(S_1,\omega_1)=\mathrm{Vol}\left(T,\frac{\omega_T}{1+r_1^2}\right)\mathrm{Vol}\left(S^2_{r_1},\frac{dA}{r_1}\right)=\left(\frac{1}{1+r_1^2}\right)\left(\frac{4\pi r_1^2}{r_1}\right)=\frac{4\pi r_1}{1+r_1^2},$$
$$\mathrm{Vol}(S_2,\omega_2)=\mathrm{Vol}\left(T,\omega_T\right)\mathrm{Vol}\left(S^2_{r_2},\frac{dA}{r_2}\right)=\left( 1 \right)\left(\frac{4\pi r_2^2}{r_2}\right)=4\pi r_2,$$
and so 
\begin{equation}
\label{Eq 9}
r_2=\frac{r_1}{1+r_1^2}.
\end{equation}

On the other hand, the second homotopy group $\pi_2(T\times S^2)=\mathbb{Z}$ with generator $\phi:S^2\rightarrow \{t\}\times S^2$, where $t\in T$ is any point.  The diffeomorphism $\psi$ of $T\times S^2$ induces an isomorphism $\psi_*$ of homotopy groups, and so $\psi_*\phi=\pm \phi$. By assumption, $\psi$ is a symplectomorphism $\psi^*\omega_2=\omega_1$, thus
$$
4\pi r_1=\int_{\phi} \omega_1 = \int_{\phi} \psi^*\omega_2= \int_{\psi_*\phi} \omega_2=\int_{\pm\phi} \omega_2=\pm 4\pi r_2.
$$
But this contradicts Equation (\ref{Eq 9}).
\end{example}

\noindent{\bf Remark 4}
 Let $\Lambda$ be a Poisson structure horizontally nondegenerate
  defined on the total space $E$ of a vector bundle $p: E \rightarrow S$
 over a compact base manifold $S$. We identify $S$ with the zero section 
and suppose that $S$ is a symplectic leaf of $(E, \Lambda)$. 
 We denote by $(\Gamma, \nu, \varphi)$ the geometric data
 associated with the first approximation $\pi$ of $\Lambda$.
 In \cite{brahic 2}, the author shows that the germ of $\Lambda$ along
 $S$ is isomorphic to a Poisson bivector field $\Lambda'$ admitting
 $(\Gamma, \nu, \varphi')$ as  associated geometric data.  The above example shows that  
$\varphi' = \varphi$ cannot occur in that case.  

\section{Acknowledgements}
Rui Loja Fernandes suggested a simplification to the proof of Example 3. Thanks also to Alan Weinstein for comments on an earlier draft of this paper.

\noindent 
 Benjamin Lent Davis\\
 Department of Mathematics and Computer Science\\
 Saint Mary's College of California\\
1928 Saint Mary's Road\\
 Moraga, CA 94556 \\
email: bldavis@stmarys-ca.edu\\[.2cm]

   A\"{\i}ssa Wade \\ 
 Department of Mathematics\\
 The Pennsylvania State University \\
 University Park, PA 16802 \\
    email: wade@math.psu.edu


\label{lastpage}
\end{document}